\documentclass[a4paper, reqno,12pt]{amsart}

\usepackage{tikz, tikz-3dplot}
\usepackage{amsmath, amsthm, amssymb,graphics }

\theoremstyle{plain}
\newtheorem{te}{Theorem}[section]

\newtheorem{co}[te]{Corollary}

\newtheorem{de}[te]{Definition}

\theoremstyle{remark}

\newtheorem*{ack*}{Acknowledgment}

\def\n{{\bf n}}

\def\0{{\bf 0}}

\def\T{{\mathbb T}}
\def\R{{\mathbb R}}

\def\C{{\mathbb C}}

\def\nint{\mathop{\diagup\kern-13.0pt\int}}

\def\Tc{{\mathcal T}}

\def\Qc{{\mathcal Q}}

\begin{document}
\author{Ciprian Demeter}
\address{Department of Mathematics, Indiana University,  Bloomington IN}
\email{demeterc@indiana.edu}

\thanks{The author is  partially supported by the  Research NSF grant DMS-1800305}
\begin{abstract} We present a slightly simpler proof of the multilinear refined Strichartz estimate from \cite{DGL}, and prove a slightly more general linear refined Strichartz estimate. Our arguments seek to clarify the connection between these estimates, refined decoupling and tube incidences.

\end{abstract}
\title[On the refined Strichartz estimates]{On the refined Strichartz estimates}

\maketitle

\section{Introduction}
Let $E$ be the extension operator associated with the paraboloid. More precisely, if $f:[-1,1]^{n-1}\to\C$, let for $x=(x_1,\ldots,x_n)$
$$Ef(x)=\int f(\xi_1,\ldots,\xi_{n-1})e(\xi_1x_1+\ldots+\xi_{n-1}x_{n-1}+(\xi_1^2+\ldots+\xi_{n-1}^2)x_n)d\xi_1\ldots d\xi_{n-1}.$$

Let us first recall the classical Strichartz estimate.

\begin{te}[\cite{Str}]\label{clasic} For each $f:[-1,1]^{n-1}\to\C$ we have
$$\|Ef\|_{L^{\frac{2(n+1)}{n-1}}(\R^{n})}\lesssim \|f\|_2.$$
\end{te}
The following result (in fact, its equivalent formulation in Corollary \ref{c1} below) was proved in \cite{DGL} when $n=2$. The proof  extends to $n\ge 3$, as observed in \cite{DGLZ}. This result proved instrumental in the recent resolution of Carleson's problem on pointwise convergence of the solution of Schr\"odinger's equation to the initial data, see \cite{DGL} and \cite{DZ}.
\medskip

The notation $\lessapprox$ will signal the existence of implicit constants of the order $O((\log R)^{O(1)})$, and $|\cdot|$ will refer to the cardinality of finite sets.

\begin{te}\label{t3}
	Let $\Omega_1,\ldots,\Omega_{n}$ be transverse cubes in $[-1,1]^{n-1}$ with diameter $\sim 1$,  in the sense that no hyperplane in $\R^{n-1}$ simultaneously intersects all $\Omega_i$. For $1\le i\le n$, let $f_i:\Omega_i\to\C$. Let $\Qc$ be a collection of $N$ pairwise disjoint cubes $q$ in $[0,R]^n$ with side length $R^{1/2}$.
	
	Then there is a subcollection $\Qc'\subset \Qc$ such that, writing $S=\bigcup_{q\in \Qc'}q$, we have
	\begin{equation}
	\label{k1}
	|\Qc|\lessapprox |\Qc'|
	\end{equation}
	and for each $\epsilon>0$
	\begin{equation}
	\label{k2}
	(\prod_{i=1}^n\|Ef_i\|_{L^{\frac{2(n+1)}{n-1}}(S)})^{1/n}\lesssim_\epsilon N^{-\frac{n-1}{n(n+1)}}R^{\epsilon}(\prod_{i=1}^n\|f_i\|_2)^{1/n}.
	\end{equation}
\end{te}

Note that the exponent $\frac{2(n+1)}{n-1}$ is the same in both theorems. The gain $N^{-\frac{n-1}{n(n+1)}}$ in Theorem \ref{t3} comes at the cost of replacing the domain of integration $\R^n$ with with a restricted collection of cubes.

We say that  a variable quantity is essentially constant if its value is always in some interval $[v,2v]$, for some fixed $v>0$. Theorem \ref{t3} admits the following corollary, which on closer inspection may in fact be seen to be an equivalent reformulation of the former, using pigeonholing as explained later in this note.
\begin{co}[Multilinear refined Strichartz estimate]
\label{c1}	
	Let $\Omega_1,\ldots,\Omega_{n}$ be transverse cubes in $[-1,1]^{n-1}$, with diameter $\sim 1$. For $1\le i\le n$, let $f_i:\Omega_i\to\C$. Let $\Qc$ be a collection of $N$ pairwise disjoint cubes $q$ in $[0,R]^n$ with side length $R^{1/2}$, such that for each $1\le i\le n$, the quantity $\|Ef_i\|_{L^{\frac{2(n+1)}{n-1}}(q)}$ is essentially constant in $q$. Then
	$$
	(\prod_{i=1}^n\|Ef_i\|_{L^{\frac{2(n+1)}{n-1}}(\cup_{q\in\Qc}q)})^{1/n}\lesssim_\epsilon N^{-\frac{n-1}{n(n+1)}}R^{\epsilon}(\prod_{i=1}^n\|f_i\|_2)^{1/n}.
	$$
\end{co}
\medskip

The proofs in \cite{DGL} and \cite{DGLZ} of this result are as follows. First, a linear refined Strichartz estimate is obtained for families of cubes $\Qc$ that have a certain structure inside $[0,R]^n$. Their cardinality is not important, only the distribution in horizontal strips plays a role. See Theorem
\ref{t5} here. Then Corollary \ref{c1} is proved by repeating this argument in the multilinear setting, pigeonholing the collection of $N$ cubes to have the additional structure from the linear setup, and using the multilinear Kakeya inequality to count ``heavy" cubes.
\medskip

In Section \ref{s2} we present a slightly different argument for Theorem \ref{t3}, one that does not rely on  a linear refined Strichartz estimate. Instead, it will use the  refined decoupling from \cite{GIOW}, that we recall below. The proof of the refined decoupling from \cite{GIOW} uses very similar ideas to the proof of the linear Strichartz estimate from \cite{DGL} and \cite{DGLZ} (in particular rescaling from $R$ to $R^{1/2}$ and the $l^2$ decoupling from \cite{BD}), but it is conceptually easier, in the sense that it relies less on creating structure, and consequently, it uses less pigeonholing.
\medskip

In Section \ref{s4} we repeat this argument to prove a slightly more general linear refined Strichartz estimate. For reader's convenience, in the last section we recall (in the simplest case $n=2$) the way the bilinear refined Strichartz estimate solves the Carleson problem.

\medskip

A tube $T$ is a cylinder in $\R^{n}$ with radius $R^{1/2}$ and length $R$. We write $R^{\epsilon}T$ to denote the cylinder with radius $R^{\frac12+\epsilon}$ and length $R$, centered at the same point as, and having the same direction as $T$. We will say that the tubes $T_1,\ldots,T_n$ with directions specified by unit vectors $\n_1,\ldots,\n_n$ are tranverse if the volume of the parallellepiped determined by these vectors has volume $\gtrsim_n1$.
\medskip

We recall the wave packet decomposition for $Ef$ on $[0,R]^n$, where $f:[-1,1]^{n-1}\to\C$
$$
Ef(x)=\sum_{T\in \Tc_R(f)}F_{T}(x)+O(R^{-100n})\|f\|_2,\;\;\;x\in [0,R]^n.
$$
The tubes in the collection $\Tc_R(f)$ intersect $[0,R]^n$.
Here $F_T=(Ef_T)1_{[0,R]^n}$, and $|F_T|$ is approximately equal to $w_T1_T$. In particular, $\|F_T\|_p\sim w_TR^{\frac{n+1}{2p}}.$ We call $w_T$ the weight of $F_T$. Moreover, the functions $f_T$ are almost orthogonal and
$$\|f\|_2\sim (\sum_{T\in\Tc_R(f)}\|f_T\|_2^2)^{1/2}.$$
See for example Chapter 2 in \cite{Dembook}.
\medskip

We now recall the two key results that will be used in this note.
\begin{te}[Refined decoupling, \cite{GIOW}]
	\label{t1}

	Let $\Qc$ be a collection of pairwise disjoint cubes $q$ in $[0,R]^n$ with side length $R^{1/2}$. Let $\epsilon>0$.
	Assume that each $q$ intersects at most $M$ fat tubes $R^\epsilon T$, with $T\in\Tc_R(f)$, for some $1\le M\lesssim_\epsilon R^{\frac{n-1}{2}+\epsilon}$.
	
	Then for each $2\le p\le \frac{2(n+1)}{n-1}$ we have	
	
	$$\|Ef\|_{L^p(\cup_{q\in\Qc}q)} \lesssim_{\epsilon}R^{O(\epsilon)}M^{\frac12-\frac1p}(\sum_{T\in\Tc_R(f)}\|F_T\|^p_{L^p(\R^n)})^{\frac1p}.$$
\end{te}
This result is a refinement of the decoupling proved in \cite{BD}, which has the factor $R^{\frac{n-1}{2}(\frac12-\frac1p)}$ in place of $M^{\frac12-\frac1p}$.

\medskip

We will also use the following multilinear Kakeya inequality from \cite{BCT}, in essentially the same way that has been used in \cite{DGL} and \cite{DGLZ}.

\begin{te}[Multilinear Kakeya]	\label{t2}
Let $\Tc_1,\ldots,\Tc_n$ be transverse families of tubes in $\R^n$. This means that $T_1,\ldots,T_n$ are transverse whenever $T_i\in\T_i$.

Given $M_1,\ldots,M_n\ge 1$, let $\Qc'$ be a collection of pairwise disjoint $R^{1/2}$-cubes $q$, with each $q$ intersecting $\sim M_i$ tubes from $\Tc_i$, for each $1\le i\le n$. Then
\begin{equation}
\label{jiefruiuir}
|\Qc'|\lesssim_\epsilon R^\epsilon\left(\frac{\prod_{i=1}^n|\Tc_i|}{\prod_{i=1}^n M_i}\right)^{\frac{1}{n-1}}.
\end{equation}
	
\end{te}	
\section{Proof that Theorem \ref{t1} and Theorem \ref{t2} imply Theorem \ref{t3} }
\label{s2}
\bigskip

We may assume that for each $q\in \Qc$ and each $1\le i\le n$ we have
\begin{equation}
\label{k3}
R^{-10n}\|f_i\|_2\le \|Ef_i\|_{L^{\frac{2(n+1)}{n-1}}(q)}\;\;,
\end{equation}
as the cubes $q$ failing to satisfy  this requirement will produce a very small contribution to the left hand side of \eqref{k2}, so they can be  harmlessly added to $\Qc'$.
\medskip

We use the wave packet decomposition for each $f_i$ on $[0,R]^n$

$$Ef_i(x)=\sum_{T_i\in \T_i}F_{T_i}(x)+O(R^{-100n})\|f_i\|_2,\;\;\;x\in [0,R]^n,$$
with $F_{T_i}=Ef_{T_i}$.
We denote by $w_{T_i}$ the weight of $F_{T_i}$.  We can partition $\T_i$ into $\lessapprox 1$  collections $\T_{i,l}$ such that $w_{T_i}$ is essentially constant for $T_i\in\T_{i,j}$ and such that
\begin{equation}
\label{k4}
Ef_i(x)=\sum_{j\lessapprox 1}\sum_{T_i\in \T_{i,j}}F_{T_i}(x)+O(R^{-100n})\|f_i\|_2,\;\;\;x\in [0,R]^n.
\end{equation}
Because of \eqref{k3} and \eqref{k4}, for each $q\in \Qc$ and each $1\le i\le n$ there is some $j_i\lessapprox 1$  such that

\begin{equation}
\label{k7}
\|Ef_i\|_{L^{\frac{2(n+1)}{n-1}}(q)}\lessapprox \|\sum_{T\in \T_{i,j_i}}F_{T_i}\|_{L^{\frac{2(n+1)}{n-1}}(q)}.
\end{equation}
We partition $\Qc$ into collections in  such a way that all $q$ in each given collection  are assigned  the same $n$-tuple $(j_1,\ldots,j_n)$. Since there are $\lessapprox 1$ such collections, we may use pigeonholing to pick one, call it $\Qc^*$, such that
\begin{equation}
\label{k5}
|\Qc|\lessapprox |\Qc^*|.
\end{equation}
Call $\T_{i,j_i}=\Tc_i$, with $(j_1,\ldots,j_n)$ being the $n$-tuple associated with $\Qc^*$, and assume the weight of each $F_{T_i}$ with $T_i\in\Tc_i$ is $\sim w_i$.
We perform one last dyadic pigeonholing.
\medskip

Fix $\epsilon>0$. We distinguish the cubes $q$ in  $\Qc^*$ according to how many fat tubes $R^{\epsilon}T$ with $T\in \T_{i,j_i}$ intersect  $q$. We may thus find a collection $\Qc'\subset \Qc^*$ such that
\begin{equation}
\label{k6}
|\Qc^*|\lessapprox |\Qc'|
\end{equation}
and such that for each $q\in \Qc'$ and each $1\le i\le n$ we have
$$M_i\le |\{T_i\in \Tc_i:\: R^\epsilon q\cap T_i\not=\emptyset\}|\le 2M_i$$
for some dyadic $M_i\lesssim R^{\frac{n-1}{2}}$. Note that \eqref{k3} combined with the Schwartz-type decay of $F_{T_i}$ away from  $T_i$ shows that $M_i$ cannot be zero.
\medskip

Using Theorem \ref{t2} (and \eqref{k5}, \eqref{k6}) we get
\begin{equation}
\label{k8}
N\lessapprox |\Qc'|\lesssim_\epsilon R^{O(\epsilon)}\left(\frac{\prod_{i=1}^n|\Tc_i|}{\prod_{i=1}^n M_i}\right)^{\frac{1}{n-1}}.
\end{equation}

Let $S=\bigcup_{q\in \Qc'}q$. Note that \eqref{k7} implies that
$$(\prod_{i=1}^n\|Ef_i\|_{L^{\frac{2(n+1)}{n-1}}(S)})^{1/n}\lessapprox (\prod_{i=1}^n\|\sum_{T_i\in\Tc_{i}}F_{T_i}\|_{L^{\frac{2(n+1)}{n-1}}(S)})^{1/n}.$$
Using Theorem \ref{t1} we can dominate the right hand side by
$$R^{\epsilon}(\prod_{i=1}^n M_i)^{\frac{1}{n(n+1)}}(\prod_{i=1}^n\sum_{T_i\in\Tc_{i}}\|F_{T_i}\|^{\frac{2(n+1)}{n-1}}_{L^{\frac{2(n+1)}{n-1}}(\R^n)})^{\frac{n-1}{2n(n+1)}}.$$

Recall that for each $p\ge 1$
$$\|F_{T_i}\|_{p}\sim w_i|T_i|^{1/p}\sim w_iR^{\frac{n+1}{2p}},$$
so the above expression is
$$\sim R^{\frac{n-1}{4}+\epsilon}(\prod_{i=1}^n M_i)^{\frac{1}{n(n+1)}}(\prod_{i=1}^n w_i)^{\frac{1}{n}}(\prod_{i=1}^n|\Tc_i|)^{\frac{n-1}{2n(n+1)}}.$$
It remains to show that
$$R^{\frac{n-1}{4}}(\prod_{i=1}^n M_i)^{\frac{1}{n(n+1)}}(\prod_{i=1}^n w_i)^{\frac{1}{n}}(\prod_{i=1}^n|\Tc_i|)^{\frac{n-1}{2n(n+1)}}\lesssim_\epsilon R^{O(\epsilon)} N^{-\frac{n-1}{n(n+1)}}(\prod_{i=1}^n\|f_i\|_2)^{1/n}.$$
Using orthogonality we write
$$\|f_i\|_{2}\gtrsim \|\sum_{T_i\in\Tc_i}f_{T_i}\|_{2}\sim w_{i}|\Tc_i|^{1/2}R^{\frac{n-1}{4}}.$$
It thus suffices to show that
$$(\prod_{i=1}^n M_i)^{\frac{1}{n(n+1)}}(\prod_{i=1}^n|\Tc_i|)^{\frac{n-1}{2n(n+1)}}\lesssim_\epsilon R^{O(\epsilon)} N^{-\frac{n-1}{n(n+1)}}(\prod_{i=1}^n|\Tc_i|)^{\frac1{2n}}.$$
After rearranging it, this is the same as \eqref{k8}.

\section{Reversing the implication?}
One may wonder if the multilinear refined Strichartz estimate (Corollary \ref{c1}) implies the refined decoupling in Theorem \ref{t1}, or at least some multilinear version of it. What we ask for is a more or less direct argument, like the one for the reverse implication described in the previous section. The answer appears to be ``no". To illustrate the relative strength of the latter compared to the former result, we point out below that the implication under question does hold if the collection $\Qc$ satisfies the saturation condition (S2) (see the statement of the following theorem). In light of \eqref{jiefruiuir}, this condition reads (we use $\approx$ to hide arbitrarily small $R^\epsilon$ factors)
$$\left(\frac{\prod_{i=1}^n|\Tc_i|}{\prod_{i=1}^n M_i}\right)^{\frac{1}{n-1}}\approx|\Qc|,$$
and means that a significant fraction of the transverse incidences between tubes occur at the cubes from $\Qc$.
\bigskip

Let $Q_1,\ldots,Q_n$ be transverse cubes in $[-1,1]^{n-1}$ and let
$f_i:Q_i\to\C$. Consider the wave packet decomposition on $[0,R]^n$

$$
Ef_i(x)=\sum_{T_i\in \Tc_R(f_i)}F_{T_i}(x)+O(R^{-100n})\|f_i\|_2,\;\;\;x\in [0,R]^n.
$$
A rather immediate computation shows that Corollary \ref{c1} implies the following result.
\begin{te}[Saturated multilinear refined decoupling]Let $p=\frac{2(n+1)}{n-1}$.
	Assume that for each $i$ and each $T_i\in \Tc_R(f_i)$, the weight of $F_{T_i}$ is $\sim 1$.
	Let $\Qc$ be a collection of pairwise disjoint cubes $q$ in $[0,R]^n$ with side length $R^{1/2}$. Let $\epsilon>0$.
	Assume that each $q$ intersects  $\sim M_i$ fat tubes $R^\epsilon T_i$, with $T_i\in\Tc_R(f_i)$, for each $1\le i\le n$. Assume that
	\\
	\\
	(S1):\; For each $i$,  $\|Ef_i\|_{L^p(q)}$ is essentially constant in $q$
	\\
	\\
	(S2):\; $$\left(\frac{\prod_{i=1}^n|\Tc_i|}{\prod_{i=1}^n M_i}\right)^{\frac{1}{n-1}}\lesssim_\epsilon R^\epsilon|\Qc|.$$

	Then
	$$(\prod_{i=1}^n\|Ef_i\|_{L^p(\cup_{q\in\Qc}q)})^{\frac1{np}}\lesssim_{\epsilon}R^{O(\epsilon)}(\prod_{i=1}^nM_i^{\frac12-\frac1p})^{\frac1n}(\prod_{i=1}^n\sum_{T_i\in\Tc_R(f_i)}\|F_{T_i}\|^p_{L^p(\R^n)})^{\frac1{np}}.$$
\end{te}
\section{Refined decoupling implies linear refined Strichartz estimate}
\label{s4}

Recall that the directions of the tubes arising in our wave packet decompositions  -let us call them {\em admissible} tubes- are given by the normal vectors to the paraboloid over $[-1,1]^{n-1}$. Thus, the angles between the directions of admissible tubes and the vertical axis $x_n$ are at most $C_n\pi$, with $C_n<\frac12$. In other words, admissible tubes are never close to being horizontal.

\begin{de}
Let $D_n$ be a fixed parameter depending only on the dimension $n$.	
A collection $\Qc$ of pairwise disjoint $R^{1/2}$-cubes in $[0,R]^n$ is said to be almost horizontal if each admissible tube intersects at most $D_n$ cubes in $\Qc$.
\end{de}

The collection of cubes intersecting a horizontal hyperplane $x_n=const$ is almost horizontal. But so is the collection of cubes intersecting the graph of a smooth function $\phi:[0,R]^{n-1}\to[0,R]$ with $\|\nabla \phi\|_{L^\infty}\lesssim_{C_n,D_n}1$.
\medskip

We prove the following result, using the refined decoupling in Theorem \ref{t1}. The reader should compare this with Theorem \ref{clasic}.

\begin{te}[Linear refined Strichartz estimate]
	\label{t5}	
	Let  $f:[-1,1]^n\to\C$. Let $\Qc$ be a collection of pairwise disjoint $R^{1/2}-$cubes $q$ in $[0,R]^n$. Assume  the quantity $\|Ef\|_{L^{\frac{2(n+1)}{n-1}}(q)}$ is essentially constant in $q$. We partition $\Qc$ into almost horizontal collections $\Qc_j$. Let $\sigma$ be the cardinality of the smallest among these collections.
	
	  Then
	\begin{equation}
	\label{dpu89erf90ew-o0-eo0f}
	\|Ef\|_{L^{\frac{2(n+1)}{n-1}}(\cup_{q\in\Qc}q)}\lesssim_\epsilon \sigma^{-\frac{1}{n+1}}R^{\epsilon}\|f\|_2.
	\end{equation}
	The implicit constant is independent of the number of collections $\Qc_j$.
\end{te}
The case of this theorem when each of the almost horizontal collections is in fact  horizontal (the cubes touch a fixed horizontal hyperplane)  was proved in \cite{DGL}. Let us briefly sketch an argument here that relies solely on Theorem \ref{t1}.
\medskip

\begin{proof}
The proof is very similar to the one in Section \ref{s2}, but it uses a trivial incidence bound, rather than the multilinear Kakeya estimate. Using pigeonholing we may assume that all the wave packets of $Ef$ have weight $\sim 1$, and that each $q\in \Qc'$ (for some $\Qc'\subset \Qc$ with $|\Qc|\lessapprox |\Qc'|$) intersects roughly $M$ fat tubes $R^\epsilon T$ with $T\in \Tc_R(f)$, for some fixed $M\ge 1$.

Then the essentially constant property combined with refined decoupling implies that
$$\|Ef\|_{L^{\frac{2(n+1)}{n-1}}(\cup_{q\in\Qc}q)}\lessapprox \|Ef\|_{L^{\frac{2(n+1)}{n-1}}(\cup_{q\in\Qc'}q)}\lesssim_\epsilon R^{\frac{n-1}{4}+\epsilon}M^{\frac1{n+1}}|\Tc_R(f)|^{\frac{n-1}{2(n+1)}}.$$
Comparing this with \eqref{dpu89erf90ew-o0-eo0f} and using that $\|f\|_2\sim R^{\frac{n-1}{4}}|\Tc_R(f)|^{1/2}$, it remains  to prove that
\begin{equation}
\label{wwdiherufuoep}
M\sigma\lessapprox |\Tc_R(f)|.
\end{equation}
It is immediate that $\Qc'$ contains an almost horizontal collection of size $\gtrapprox \sigma$, call it $\Qc^*$. Indeed, since
$$\sum_j|\Qc'\cap \Qc_j|=|\Qc'|\gtrapprox |\Qc|=\sum_j|\Qc\cap \Qc_j|$$
there must be some $j$ with $|\Qc'\cap \Qc_j|\gtrapprox |\Qc\cap \Qc_j|\ge \sigma$. We let $\Qc^*=\Qc'\cap \Qc_j$.

Recall that each $q\in \Qc^*$ intersects roughly $M$ tubes, and that each tube can intersect at most $D_n$ cubes in $\Qc^*$. It follows that
$$M|\Qc^*|\lesssim  D_n|\Tc_R(f)|.$$

The estimate \eqref{wwdiherufuoep} is now immediate.

\end{proof}
\section{Carleson's problem}
\label{s5}
Carleson's problem \cite{Car} in the plane about the pointwise convergence of the solution to Schr\"odinger's equation to initial data is equivalent to the following theorem.

\begin{te}
Let $f_1,f_2$ be supported on disjoint intervals in $[-1,1]$. Consider a collection of $\sim R$ unit squares $\omega$ in $[0,R]^2$ with pairwise disjoint projections onto the $x$ axis, and let $S$ be their union. Then
$$\|(Ef_1Ef_2)^{1/2}\|_{L^2(S)}\lesssim R^{\frac14+\epsilon}(\|f_1\|_2\|f_2\|_2)^{\frac12}.$$
\end{te}
Let us see how Theorem \ref{t3} implies this. Pigeonholing, we may assume that
$\|(Ef_1Ef_2)^{1/2}\|_{L^2(\omega)}$ is essentially constant in $\omega$, and that there are either 0 or $\sim \lambda$ squares $\omega$ inside each $R^{1/2}$-square $q$ (we take $q$ from a partition of $[0,R]^2$). The value of $\lambda$ will be irrelevant. Let $\Qc$ be the collection of $R^{1/2}$-squares containing the squares $\omega$. Assume it has size $N$. Apply Theorem \ref{t3} to find $\Qc'\subset \Qc$ such that $N\lessapprox |\Qc'|$ and
$$\|(Ef_1Ef_2)^{1/2}\|_{L^6(\cup_{q\in \Qc'}q)}\lesssim_\epsilon R^{\epsilon}N^{-1/6}(\|f_1\|_2\|f_2\|_2)^{\frac12}.$$
Note that because of the essentially constant assumption and the $\lambda$ uniformity we have
$$\|(Ef_1Ef_2)^{1/2}\|_{L^2(S)}\lessapprox \|(Ef_1Ef_2)^{1/2}\|_{L^2(\cup_{q\in \Qc'}\cup_{\omega\subset q}\omega)}.$$
Let us denote by $M$ the cardinality  of $\{\omega\subset q:\;q\in\Qc'\}$. Recall that $M\le R$. Using H\"older, the right hand side is dominated by
$$M^{1/3}\|(Ef_1Ef_2)^{1/2}\|_{L^6(\cup_{q\in \Qc'}\cup_{\omega\subset q}\omega)}.$$
This is trivially dominated by
$$M^{1/3}\|(Ef_1Ef_2)^{1/2}\|_{L^6(\cup_{q\in\Qc'}q)}.$$
Note also that
$$M\approx \lambda N,\;\;\text{so}\;\; N\gtrapprox MR^{-1/2}.$$
Combining these leads to the desired estimate
$$\|(Ef_1Ef_2)^{1/2}\|_{L^2(S)}\lesssim_\epsilon R^{\epsilon }M^{1/3}N^{-1/6}(\|f_1\|_2\|f_2\|_2)^{\frac12}\lesssim R^{\frac14+\epsilon}(\|f_1\|_2\|f_2\|_2)^{\frac12}.$$


\begin{thebibliography}{99}
	\bibitem{BCT} Bennett, J.,  Carbery, A. and  Tao, T. {\em On the multilinear restriction and Kakeya conjectures}, Acta Math. 196 (2006), no. 2, 261-302
	\bibitem{BD}Bourgain, J. and Demeter, C. {\em The proof of the $l^2$ Decoupling Conjecture}, Annals of Math. 182 (2015), no. 1, 351-389.
		\bibitem{Dembook} Demeter, C. {\em Fourier restriction, decoupling and applications},  Cambridge University Press, 2020
	\bibitem{Car} Carleson, L. {\em Some analytic problems related to statistical mechanics} In: Euclidean harmonic analysis (Proc. Sem., Univ. Maryland, College Park, Md., 1979), pp. 5-45. Lecture Notes in Math., vol. 779. Springer, Berlin 	
	\bibitem{DZ} Du, X and Zhang R {\em Sharp $L^2$ estimate of Schr\"odinger maximal function in higher dimensions},  Ann. of Math. (2) 189 (2019), no. 3, 837-861
	\bibitem{DGL} Du, X, Guth, L and Li, X {\em A sharp Schr\"odinger maximal estimate in $\R^2$} Ann. of Math. (2) 186 (2017), no. 2, 607-640
	\bibitem{DGLZ} Du, X, Guth, L,  Li, X and Zhang, R {\em Pointwise convergence of Schr\"odinger solutions and multilinear refined Strichartz estimates},  Forum Math. Sigma 6 (2018), e14, 18 pp
	\bibitem{GIOW} Guth, L.,  Iosevich, A., Ou, Y. and Wang, H.  {\em On Falconer's distance set problem in the plane}, to appear in Invent. Math.
	\bibitem{Str} Strichartz, R. S. {\em Restrictions of Fourier transforms to quadratic surfaces and decay of solutions of wave equations}, Duke Math. J. 44 (1977), no. 3, 705-714.
\end{thebibliography}
\end{document}